\documentclass[pdflatex,sn-mathphys-num]{sn-jnl}

\usepackage{tikz}
\usepackage{pifont} 
\usepackage{xypic}
\usepackage{graphicx}%
\usepackage{multirow}%
\usepackage{amsmath,amssymb,amsfonts}%
\usepackage{amsthm}%
\usepackage{mathrsfs}%
\usepackage[title]{appendix}%
\usepackage{xcolor}%
\usepackage{textcomp}%
\usepackage{manyfoot}%
\usepackage{booktabs}%
\usepackage{algorithm}%
\usepackage{algorithmicx}%
\usepackage{algpseudocode}%
\usepackage{soul}
\usepackage{listings}%
\usepackage[font=small,skip=0pt]{caption}
\usepackage[font=small,skip=0pt]{subcaption}
\usepackage{tabularx}
\usepackage{longtable}
\usetikzlibrary{arrows.meta, positioning}
\usetikzlibrary{shapes.geometric}

\theoremstyle{thmstyleone}%
\theoremstyle{remark}

\newtheorem{theme}{Research Theme}
\newtheorem{contribution}{Contribution}
\newtheorem{step}{Documentation Step}

\theoremstyle{thmstyletwo}%

\theoremstyle{thmstylethree}%

\begin{document}

\title[]{Counting Who Counts: Black Women Mathematicians in the United States}


\author[1]{\fnm{Olayemi} \sur{Adeyemi}}\email{olayemiadeyemi@u.boisestate.edu}
\equalcont{Authors are listed in alphabetical order.}

\author[2]{\fnm{Ron} \sur{Buckmire}}\email{ron.buckmire@marist.edu}
\equalcont{Authors are listed in alphabetical order.}

\author*[3]{\fnm{Ranthony A.} \sur{ Clark}}\email{edmonds.110@osu.edu}
\equalcont{Authors are listed in alphabetical order.}

\author[4]{\fnm{Shanzhen} \sur{Gao}}\email{sgao@vsu.edu}
\equalcont{Authors are listed in alphabetical order.}

\author[5]{\fnm{Jakini} \sur{Kauba}}\email{jakauba@gmail.com}
\equalcont{Authors are listed in alphabetical order.}

\author[6]{\fnm{Lakeshia} \sur{Legette-Jones}}\email{ljones3@cau.edu}
\equalcont{Authors are listed in alphabetical order.}

\author[7]{\fnm{Apoorva} \sur{Mate}}\email{amate@psu.edu}
\equalcont{Authors are listed in alphabetical order.}

\author[8]{\fnm{Omayra} \sur{Ortega}}\email{ortegao@sonoma.edu}
\equalcont{Authors are listed in alphabetical order.}

\author[9]{\fnm{José} \sur{Pabón}}\email{jlpabon@alumni.princeton.edu}
\equalcont{Authors are listed in alphabetical order.}

\author[10]{\fnm{Robin T.} \sur{Wilson}}\email{robin.wilson@lmu.edu}
 \equalcont{Authors are listed in alphabetical order.}

\affil[1]{\orgdiv{Department of Mathematics}, \orgname{Boise State University}, \orgaddress{\city{Boise}, \postcode{83725}, \state{ID}, \country{USA}}}

\affil[2]{\orgdiv{School of Computer Science and Mathematics}, \orgname{Marist University}, \orgaddress{\city{Poughkeepsie}, \postcode{12601}, \state{NY}, \country{USA}}}

\affil*[3]{\orgdiv{Department of Mathematics}, \orgname{The Ohio State University}, \orgaddress{\city{Columbus}, \postcode{43210}, \state{OH}, \country{USA}}}

\affil[5]{\orgdiv{Reginald F. Lewis College of Business}, \orgname{Virginia State University}, \orgaddress{\city{Petersburg}, \postcode{23806}, \state{VA}, \country{USA}}}

\affil[5]{\orgdiv{School of Mathematical and Statistical Sciences}, \orgname{Clemson University}, \orgaddress{\city{Clemson}, \postcode{29634}, \state{SC}, \country{USA}}}

\affil[6]{\orgdiv{Department of Mathematical Sciences}, \orgname{Clark Atlanta University}, \orgaddress{\city{Atlanta}, \postcode{30314}, \state{GA}, \country{USA}}}

\affil[7]{\orgname{Pennsylvania State University -- Brandywine Campus}, \orgaddress{\city{Media}}, \postcode{19063}, \state{PA}, \country{USA}}

\affil[8]{\orgdiv{Mathematics and Statistics Department}, \orgname{Sonoma State University}, \orgaddress{\city{Rohnert Park}}, \postcode{94928}, \state{CA}, \country{USA}}

 \affil[9]{\orgdiv{Mathematics Department}, \orgname{Princeton University}, \orgaddress{\city{Princeton}, \postcode{08540}, \state{NJ}, \country{USA}}}

 \affil[10]{\orgdiv{Seaver College of Science and Engineering}, \orgname{Loyola Marymount University}, \orgaddress{\city{Los Angeles}}, \postcode{90045}, \state{CA}, \country{USA}}

\abstract{For Black women in the mathematical sciences, histories have frequently been fragmented or omitted from official records, an absence that itself constitutes a form of exclusion. This project addresses this critical gap by constructing a curated dataset of 147 Black women who have earned doctoral degrees in mathematics and related mathematical sciences in the United States. By integrating archival research, institutional records, digital databases, and community knowledge, we establish a reproducible methodology for identifying, verifying, and documenting Black women mathematicians.\\ 

We situate this dataset within a broader history of efforts to document Black mathematicians and describe how prior archival and community based initiatives informed its development. We then demonstrate the utility of the dataset through exploratory analyses of doctoral production, institutional pathways, undergraduate preparation, mentorship networks, and research areas. This revealed patterns including large concentrations of doctoral production among a small number of institutions and advisors.
Additionally, our analysis highlights the vital role of Historically Black Colleges and Universities (HBCUs) in preparing Black women who have earned doctoral degrees in the mathematical sciences. \\

By transforming historical documentation into structured research infrastructure, this work provides a foundation for future investigations of mathematical communities using historical as well as mathematical, statistical, and computational methods. More broadly, the framework developed here offers a model for constructing datasets that support the study of underrepresented communities within the mathematical sciences.
}

\keywords{academic genealogy; Mathematical community; history of mathematics, data curation; archival research; institutional pathways; doctoral education; representation in mathematics; Black women}

\maketitle

\section{Introduction}

Understanding a mathematical community begins with knowing who belongs to it. While mathematics is often viewed as a discipline built on abstraction and universal principles, the people who create, teach, and advance mathematical knowledge exist within social and historical contexts that shape their opportunities and experiences. For Black women in the mathematical sciences, those histories have frequently been fragmented, under documented, or omitted from the official records of the discipline. The absence of data is itself a form of exclusion. When a population is not systematically recorded, it becomes difficult to study, support, fund, or even recognize its contributions. While there have been historical efforts to document Black mathematicians and Black women in mathematics in particular, challenges related to data validation, availability, and research infrastructure have prevented the creation of a complete, not partial, Census of Black women in mathematics.  

This project addresses that gap by constructing a curated dataset of Black women who have earned doctoral degrees in mathematics and related mathematical sciences in the United States. The significance of this work extends beyond documentation. By creating a verified historical record, this project preserves the contributions of scholars whose collective impact cannot be fully understood through individual biographies alone, and lays the foundation for future studies on patterns of matriculation in mathematics related to the community of Black women mathematicians.

\subsection{Why Count?}

The mathematical community cannot effectively support, mentor, or advocate for its members without first knowing who they are and how they are connected. Documentation can shape how individuals understand the mathematical sciences. Students, educators, researchers, and the broader public often form their sense of belonging in mathematics and their understanding of the broader mathematical community through the people whose stories they encounter. An incomplete historical record can obscure contributions of entire communities and unintentionally reinforce narrow perceptions of who can be a mathematician. By documenting Black women mathematicians, this work makes visible their contributions and longstanding presence within the mathematical sciences.

Similarly, researchers cannot study patterns of participation, institutional pathways, or professional trajectories without systematic information about the individuals who comprise a field. Counting and documenting mathematicians establishes a record from which questions of participation, pathways, and contributions can be studied. By documenting the presence and growth of Black women in the mathematical sciences, this project provides a foundation for future community building efforts, mentorship initiatives, and evidence based policy discussions.

The importance of documenting mathematical communities is reflected in decades of efforts to preserve the histories of Black mathematicians. Foundational works such as ``Black Mathematicians and their Works" \cite{newell1980black}, Patricia Kenschaft's study of Black women mathematicians \cite{MR628027}, Erica Walker's \emph{Beyond Bannekar} \cite{walker2014beyond}, Scott Williams' Mathematicians of the African Diaspora (MAD) project \cite{madpages}, and MADDER \cite{goins2021} preserved names, biographies, and professional accomplishments that might otherwise have been overlooked. These efforts also demonstrated the difficulty of constructing a complete historical record and the need for continued documentation. Section~\ref{sec:documenting} provides a detailed discussion of these initiatives and their influence on the present work.

\subsection{Why Now?} 

The project was motivated by the centennial of Elbert Frank Cox's 1925 achievement of becoming the first Black person to receive a doctorate in mathematics and discussions emerging from the 2025 workshop on \emph{MetaMath} at the American Institute for Mathematics (AIM). MetaMath seeks to study the mathematical community using mathematical tools, and our ultimate goal is to utilize techniques from mathematics, statistics, and data science to study the community of Black women mathematicians. In order to do this we first needed to investigate, excavate, and most importantly enumerate as complete a record as possible of Black women who have received doctoral degrees in mathematics. 

The present moment represents a unique opportunity to count and document Black women in mathematics. Over the last century, the tools available for documenting Black mathematicians have changed substantially. Decades of documentation efforts (e.g. \cite{houston2025, goins2021, madpages, newell1980black}) coupled with modern archival resources make the systematic enumeration of Black mathematicians possible in ways that were previously unattainable. To date we have identified 147 Black women mathematicians, and at the time of the writing of this article we continue to update and curate our growing dataset. 

We aim to transform historical documentation from an archival record into critical research infrastructure capable of supporting mathematical, statistical, historical, and computational investigations of the mathematical sciences community. A curated dataset supports the study of long term trends in doctoral production, institutional participation, geographic distribution, advising networks, and professional pathways. These analyses provide insight into how mathematical communities develop over time and create opportunities for future quantitative and historical investigations.

\subsection{Contributions}

This work extends previous documentation efforts by developing a reproducible framework for constructing and analyzing a dataset of Black women mathematicians. The primary contributions described in this manuscript are:

\begin{contribution}
A synthesis of historical efforts to document Black women mathematicians (Section~\ref{sec:documenting}).
\end{contribution}

\begin{contribution}
A reproducible methodology for identifying, verifying, and documenting Black women mathematicians (Section~\ref{sec:methods}).
\end{contribution}

\begin{contribution}
A curated and continuously evolving dataset of Black women mathematicians in the United States, including educational, institutional, and research metadata (Section~\ref{sec:methods}). The dataset is designed to incorporate additional records and refinements as new sources and community contributions are identified.\footnote{We provide a link to a form for (self) reporting Black mathematicians: \url{https://tinyurl.com/psv6cm2e}. At the time of writing, additional calls for participation have been published in the Spring 2026 Newsletter of the National Association of Mathematicians \cite{nam2026phdsdiaspora} and the May--June newsletter of the Association for Women in Mathematics \cite{awm2026enumerating}.}
\end{contribution}

\begin{contribution}
A demonstration of how historical documentation can serve as research infrastructure for quantitative study of mathematical communities (Sections~\ref{sec:results}--\ref{sec:conclusion}).
\end{contribution}

Together, these contributions establish a foundation for future research on Black mathematicians, and provide a general framework that can be adopted to study other mathematical communities.

\subsection{Organization}

In Section~\ref{sec:documenting} we provide a loose timeline of select efforts related to documenting and highlighting Black women mathematicians in the United States that preceded and influenced the study presented in this article. Section~\ref{sec:methods} details how our data was collected, common sources, and detailed metadata about the ongoing development of our dataset of Black women mathematicians. There we also include discussions on how identity was handled, what we mean by `Black' and by `mathematician.' In Section~\ref{sec:results} we describe emerging stories from the early explorations of our growing dataset. There we focus on patterns of high concentration of graduates by institution type, advisor, and research area. Lastly, in Section~\ref{sec:conclusion} we discuss the future work for a systemic study of Black mathematicians that is enabled by our data collection and curation. In addition, we discuss plans for parallel efforts that include Black men in mathematics. We conclude with the broader context of this work with an emphasis of how utilizing mathematics to study the mathematical community ties \emph{Metamath} to the framework of the emerging research field of \emph{quantitative justice}.

\section{Documenting and Counting Black Mathematicians}\label{sec:documenting}

Here we present a truncated timeline of relevant efforts to document the production of mathematical doctoral degrees in the United States. We incorporate efforts that include underrepresented mathematicians in general, with an emphasis on those projects which specifically focus on Black women mathematicians. Together, these digital archives, research databases, and historical scholarship show how the documentation of mathematicians has evolved over the past four decades to lay the foundation for our present work.

\begin{table}[htpb]
\centering
\textbf{\normalsize Historical Initiatives to Enumerate and Document Black Mathematicians} \\ [0.75em]

\begin{tabularx}{\textwidth}{p{0.5cm} >{\raggedright\arraybackslash}p{2.3cm} >{\raggedright\arraybackslash}X p{1.3cm} >{\raggedright\arraybackslash}p{1.5cm}}
\toprule
\textbf{Year} & \textbf{Initiative} & \textbf{Contribution} & \textbf{Format} & \textbf{Counts} \\
\midrule
1980 & ``Black Mathematicians and their Works" & First national reference volume documenting Black mathematicians and mathematics educators & Reference volume & 62 mathematicians \\
\addlinespace
1981 & ``Black Women in Mathematics" & First focused documentation of Black women PhD mathematicians & Journal article & 21 Black women \\
\addlinespace
2014 & \textit{Beyond Banneker} & Historical and qualitative study of Black mathematicians' educational and professional experiences & Monograph & 35 interviews \\
\addlinespace
2020 & MADDER & Expansion of biographical documentation into a searchable archive & Digital database & 700+ biographies \\
\addlinespace
2022 & CURM Black Women in Mathematics Project & Expansion and correction of the historical record of Black women PhD mathematicians & Curated dataset & 77 Black women \\
\addlinespace
2025 & AIM MetaMath Workshop  & Synthesis of historical documentation into a dataset for quantitative investigation & Curated dataset & 119 Black women \\
\addlinespace
2026 & Documenting and Enumerating Black Women Mathematicians & Expansion of AIM MetaMath Workshop dataset including additional verified mathematicians and metadata & Curated dataset & 147 Black women \\
\bottomrule
\end{tabularx}

\caption{Selected milestones in the documentation and counting of Black mathematicians in the United States. Documentation efforts evolved from printed reference works to digital databases and curated research datasets capable of supporting quantitative analyses.}
\label{tab:documentation_timeline}

\end{table}

\subsection{Black Mathematicians and their Works (1980)}\label{sec:bmatw}

The earliest large scale effort to systematically document Black mathematicians in the United States was the 1980 volume ``Black Mathematicians and their Works" \cite{newell1980black}, edited by Virginia K. Newell, Joella H. Gipson, L. Waldo Rich, and Beauregard Stubblefield. They were explicit in their inclusion requirements: ``We are using as criteria Black mathematicians who has earned terminal degrees in pure or related fields in mathematics education or mathematicians who are outstanding in their fields." It highlighted important pioneers, noting Benjamin Banneker as the first widely recognized Black mathematician in the United States, Elbert Frank Cox as the first Black American to earn a PhD in mathematics, and Evelyn Boyd Granville and Marjorie Lee Browne as the first Black women to receive doctorates in mathematics in 1949\footnote{During this year there were only 147 PhDs granted in mathematics in the United States \cite{newell1980black}.}. In total, they identified sixty two Black mathematicians and educators whose careers spanned from 1731 (Banneker) through doctoral recipients of the early 1970s. In addition to brief biographical sketches, selected publications were included in the volume to highlight the academic achievements of Black mathematicians\footnote{Almost half of these individuals had published scholarly work, with three mathematicians, Ernest Wilkins Jr., David Blackwell, and A. T. Bharucha-Reid, authoring over forty research papers. David Blackwell was also noted as the only Black mathematicians who had been elected to the National Academy of Sciences\cite{newell1980black}.}.

\begin{figure}[htpb]
    \centering

     \textbf{Black Mathematicians and Their Works}\par\vspace{0.5em}
    
    \begin{subfigure}{0.425\textwidth}
        \centering
        \includegraphics[
            width=\linewidth,
            height=6in,
            keepaspectratio
        ]{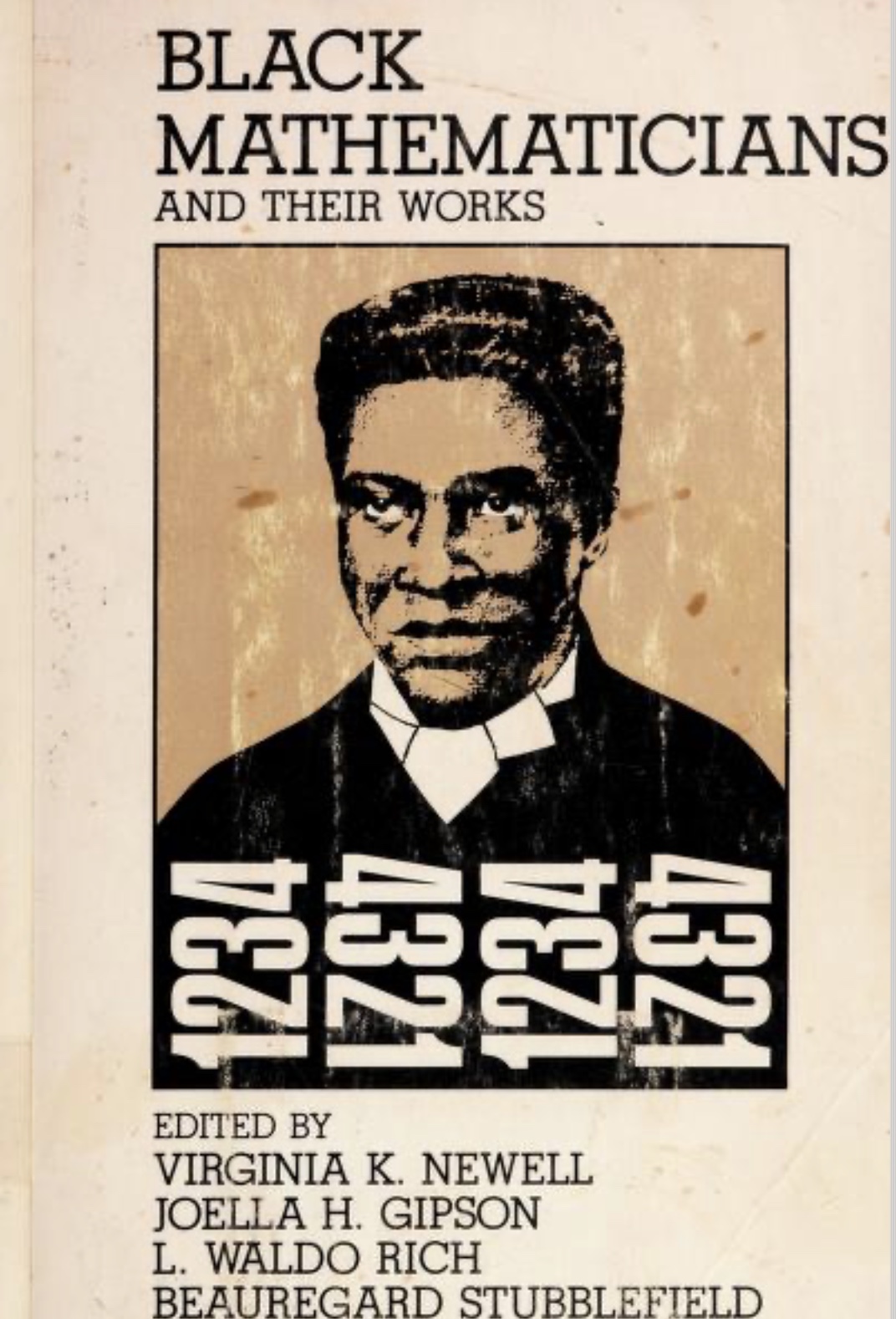}

    \end{subfigure}
    \begin{subfigure}{0.55\textwidth}
        \centering
        \includegraphics[
            width=\linewidth,
            height=6in,
            keepaspectratio
        ]{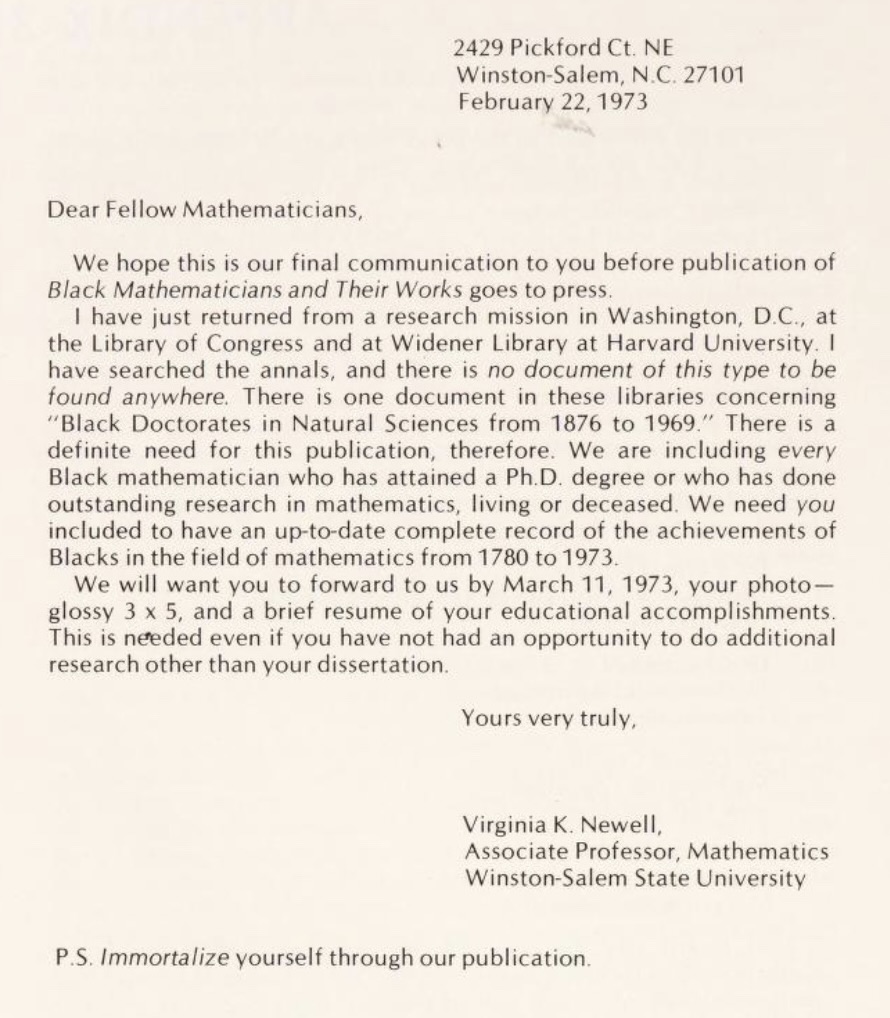}
    \end{subfigure}

    \caption{The cover of ``Black Mathematicians and their Works" \cite{newell1980black} (left) and an excerpt from an accompanying solicitation letter sent to Black mathematicians requesting photographs and biographical information for the volume (right).}
    
    \label{fig:bmatw}
\end{figure}

Of the sixty two Black mathematicians appearing in this work, thirteen were Black women\footnote{In the preface of \cite{newell1980black}, it was noted that between 1920 and 1943, 4,942 of the 171,647 (approximately $2.9\%$) doctorates awarded by universities in the United States were in mathematics. If the time period was focused to when the first eight Black PhDs in mathematics were awarded, 1875 to 1943, there were 1,374 total PhDs awarded in mathematics in the United States.} (approximately $21\%$). The editors used a list of 100 mathematicians compiled by the National Association of Mathematicians (NAM), an organization founded in 1969 whose chief aims included promoting the mathematical development and excellence of all underrepresented minorities. Individuals on this list were sent letters to solicit correspondence for inclusion in their book. These letters highlighted the editors' recognition that no comparable reference work on Black mathematicians existed, as well as their vision of creating a permanent historical record through community participation (Figure~\ref{fig:bmatw}).

\subsection{Kenschaft's List (1981)}\label{sec:kenschaft}

One of the first published works specifically related to Black women mathematicians was Patricia Kenschaft’s 1981 article in the American Mathematical Monthly titled ``Black Women in Mathematics in the United States” \cite{MR628027}.  In this article, Kenschaft shared names and brief biographical information about Black women in the United States holding doctoral degrees in mathematics at the time. She identified twenty one Black women with PhDs in mathematics\footnote{The research efforts for Kenschaft's article on Black women in math predated the publication of the volume described in Section~\ref{sec:bmatw}. Her references include conference activities by the Association for Women in Mathematics as well as the 1970 decennial census. From the census she was able to determine that 1,100 Black women self reported as mathematicians, 244 of those women reported that they were post secondary instructors, 178 had finished four or more years of college, and 152 had completed five or more years of college \cite{MR628027}.}. For some time, this list of twenty one was thought to be complete until 2001 when Euphemia Lofton Haynes, who was not featured in the article, was identified as the first Black Women to earn her PhD in mathematics from Catholic University in Washington, D.C. in 1943. 

Kenschaft continued to study Black women in mathematics \cite{kenschaft1993blackwomen} and Black mathematicians in general \cite{kenschaft1987blackresearch, kenschaft1996whatnext}. Her 1981 article however, served as an important historical benchmark that both inspired others to investigate Black women mathematicians and modeled the importance of providing biographical information in historical documentation efforts whenever possible. Four decades later, the Black Women in Mathematics CURM project (Section~\ref{sec:CURM}) utilized Kenschaft's list as the impetus for research. Given this baseline, the team considered who was omitted from the original list, and who had earned doctoral degrees after 1981. In investigating these questions, this project sought to build a more historically complete record of Black women in mathematics.

\subsection{\textit{Beyond Banneker} (2014) }\label{sec:banneker_archives}

The seminal text, \textit{Beyond Banneker: Black Mathematicians and the Paths to Excellence}, by Dr. Erica N. Walker, was published in 2014 \cite{walker2014beyond}. Unlike ``Black Mathematicians and their Works" which emphasized the mathematical contributions of Black mathematicians in general, and Kenschaft's article which offered short biographical sketches of Black women mathematicians, \emph{Beyond Banneker} studied historical factors influencing formative, educational, and professional experiences in mathematics. Walker interviewed thirty five Black mathematicians who earned their doctoral degrees between 1940 and 2010. 

Walker sorted the discussion of these mathematicians into significant time periods based on degree production--the First Generation (1940 to 1965), the Second Generation (1966 to 1985), and the Third Generation (1986 to 2010). Each period reflected different factors influencing the trajectories and experiences of Black mathematicians, with the first generation doing mathematics during de facto segregation and the dawn of the Civil Rights movement, the second navigating de jure segregation, desegregation, and affirmative action programs, and the third developing careers in the modern 21st century academic landscape. In 2017 Erica Walker published ``Excellence and devotion: black women in mathematics in the United States” \cite{MR3775375}, which extended on her previous work by focusing on qualitative interviews with Black women to study the intersection of race and gender in mathematics.

\subsection{MAD to MADDER  (2020) }\label{sec:mad}

In 1997, the Mathematicians of the African Diaspora (MAD pages) at the University at Buffalo was launched as an ambitious attempt to archive every Black person who earned a PhD in the mathematical sciences in the United States. This marked a transition from print based documentation efforts to publicly accessible digital archives of Black mathematicians. This was a labor of love by Scott W.\ Williams, who began this site during his tenure as a Professor of Mathematics at the University of Buffalo. This effort relied upon the network of mathematicians in National Association of Mathematicians organization and included individuals that Williams met personally in his travels to conferences, as well as others that were shared with him via email and written correspondence. All together, Williams put together more than seven hundred biographical profiles \citep{goins2021} of Black scientists in mathematics, computer science, and physics. After retiring in 2008, the site remained online, but was no longer updated.

In 2015, Edray Herber Goins, Donald King, Asamoah Nkwanta, and others worked to develop Mathematicians of the African Diaspora Database’s Ensemble of Researchers, MADDER, an updated version of the original MAD pages. Unlike the original MAD pages, MADDER was designed not only to preserve existing biographies but also to support ongoing data collection through a searchable, expandable research database. This updated site went public on October~9, 2020 \citep{goins2021}. MADDER preserved Williams' vision while providing infrastructure for future documentation projects of Black mathematicians, including the Black Women in Mathematics CURM initiative described below (Section~\ref{sec:CURM}).

\subsection{CURM at Pomona College and Cal Poly Pomona (2022)}\label{sec:CURM}

In 2022 a group of undergraduate students from Pomona College and California State Polytechnic University, Pomona (Cal Poly Pomona) conducted an undergraduate research project titled ``Black Women in Mathematics” under the supervision of Edray Goins and Robin Wilson. This project was funded by the Center for Undergraduate Research in Mathematics (CURM), and was part of the ongoing project MADDER, led by Goins, that focused on tracking and documenting the contributions of mathematicians of the African Diaspora (See Section~\ref{sec:mad}). The goal of the project was to construct the most complete record possible of Black women PhDs in mathematics.

The project began with Patricia Kenschaft's 1981 list as a historical baseline. The Pomona group identified thirty two Black women that received their PhD in mathematics before 1981 and at least forty seven that earned their PhD after 1981. They identified ten women in addition to Euphemia Lofton Haynes that received their PhDs before 1981 that were not listed in the Patricia Kenschaft article (Sadie Gasaway, Anne Marie Watkins Garraway, Anne Marie Meyers Pelot, Gloria Ford Gilmer, Karolyn Anne Morgan, Louise Nixon Sutton, Rita Louise Dixion, Mary C. Embry-Wardrop, Beryl E. Hunte, and Emma Fenceroy.\begin{footnote}{Note that an Errate for this publication was published in 1982 that included Emma Fenceroy and Anne Garraway \cite{MR1540102}.}\end{footnote}$^{,}$\begin{footnote}{Subsequent research revealed that the following women also earned their Ph.D.s before 1981: Louise Nixon Sutton, Rita Louise Dixion, Mary C. Embry-Wardrop, and Beryl E. Hunte.}\end{footnote} They utilized another article, ``An Unlikely Cohort: Disrupting the Stereotypes about Mathematicians and their Careers” \cite{MR4378418} which listed Black women degree recipients from the University of Maryland, and identified one person (Evelyn Patterson Scott) that was not originally included in their list of the first 20 Black women to earn doctoral degrees in mathematics in the United States. In addition, they pulled upon existing archives, databases, and institutional contacts to create an updated list of Black women mathematicians in the United States. This data was ultimately compiled into a spreadsheet of seventy seven names which became the initial database for the project described in this article.

\subsection{AIM Workshop on MetaMath  (2025)}\label{sec:AIM}

In December 2025, the American Institute for Mathematics (AIM) hosted a workshop entitled \emph{MetaMath: Modeling the mathematical sciences community using mathematics, statistics, and data science}.  Modeled after the Science of Science 
\cite{fortunato2018science}, MetaMath seeks to use  mathematics to study the mathematical sciences community itself \cite{Buckmire2025MathematicsOfMathematics}. Among the workshop participants were researchers already engaged in documenting Black mathematicians. Combined with the centennial of Elbert Frank Cox's PhD \cite{houston2025}, this created a natural opportunity to examine Black mathematicians as a case study in MetaMath. The goal was to utilize methods from mathematics, statistics, and data science to study Black mathematical communities. 

In order to move towards this, it was first necessary to synthesize historical efforts of documentation. The foundation for this research was the data product of seventy seven women from the CURM ``Black Women in Mathematics" project. The number of Black women identified was increased dramatically, with a count of one hundred nineteen by the end of the workshop. It is important to note that the workshop marked a historical shift from documentation of Black mathematicians to quantitative investigation. Given a substantially expanded and curated dataset of Black women in mathematics, it is then possible to study patterns of doctoral production, institutional pathways, geographic mobility, mentorship networks, and other characteristics of the mathematical sciences community. The data collection effort described in this article allowed our team to begin to turn an archival resource into research infrastructure capable of supporting mathematical, statistical, and historical analyses.

The transition from documentation to quantitative investigation requires a dataset that is both historically grounded and methodologically rigorous. Accordingly, the remainder of this paper focuses on the construction of the dataset itself, which includes one hundred forty seven Black women. We include the criteria used to identify individuals, the sources consulted, and the procedures used to collect, verify, and organize the data for future mathematical, statistical, and historical analyses.

\section{Data Collection and Methodology}\label{sec:methods}

The dataset was developed using the documentation framework illustrated in Figure~\ref{fig:documentation_pipeline} below. The following subsections describe the inclusion criteria, data sources, collection procedures, and validation methods used to construct the dataset.

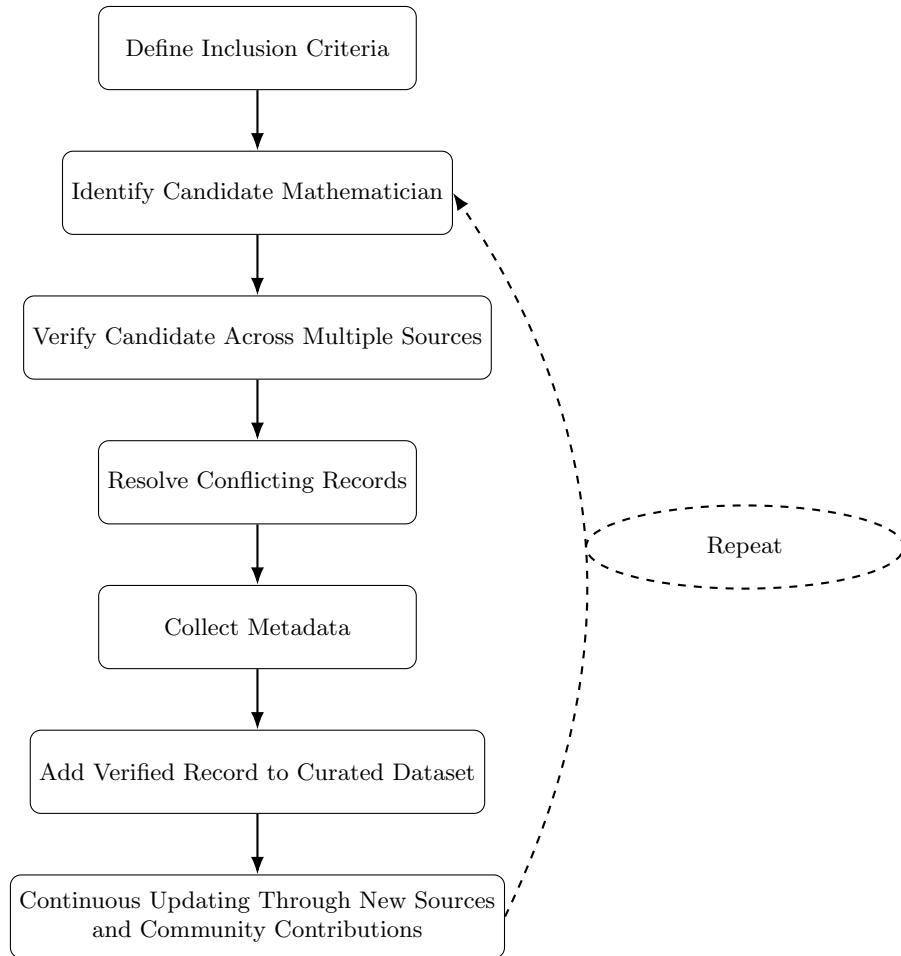
\begin{figure}[h]
\centering
\textbf{\normalsize Framework Visualization  for Documenting Black Women Mathematicians} \\ [0.75em]

\begin{tikzpicture}[
    node distance=0.8cm,
    >=Latex,
    every node/.style={
        draw,
        rectangle,
        rounded corners,
        align=center,
        minimum width=4.2cm,
        minimum height=1.1cm,
        font=\small
    },
    arrow/.style={->, thick}
]

\node (criteria) {Define Inclusion Criteria};

\node (discover) [below=of criteria]
{Identify Candidate Mathematician};

\node (verify) [below=of discover]
{Verify Candidate Across Multiple Sources};

\node (resolve) [below=of verify]
{Resolve Conflicting Records};

\node (metadata) [below=of resolve]
{Collect Metadata};

\node (dataset) [below=of metadata]
{Add Verified Record to Curated Dataset};

\node (update) [below=of dataset]
{Continuous Updating Through New Sources\\
and Community Contributions};

\draw[arrow] (criteria) -- (discover);
\draw[arrow] (discover) -- (verify);
\draw[arrow] (verify) -- (resolve);
\draw[arrow] (resolve) -- (metadata);
\draw[arrow] (metadata) -- (dataset);
\draw[arrow] (dataset) -- (update);

\draw[->, dashed, thick, bend right=30]
(update.east) to node[
    right,
    shape=ellipse,
    draw,
    fill=white,
    inner sep=1.5pt,
    font=\small
]{Repeat}
(discover.east);

\end{tikzpicture}

\caption{Documentation workflow used to construct and continuously update the curated dataset of Black women mathematicians.}
\label{fig:documentation_pipeline}

\end{figure}

\subsection{Inclusion Criteria}

This project documents Black women who earned doctoral degrees in mathematics and related mathematical sciences in the United States. Because historical records of mathematicians have not always used consistent disciplinary boundaries, institutional classifications, or racial identity categories, we clarify the criteria used to construct this dataset. Here, we focus on what is meant by ``Black" and by ``mathematician."

In this project, \emph{Black} refers to a self identified racial category that includes African Americans, Afro-Caribbean, Afro-Latino, and members of the broader continental African diaspora. The category is not restricted to a particular national origin or geographic region. Where archival records predate self reporting, we use contextual inference from sources such as institutional records, family histories, and professional organization affiliations (Section~\ref{sec:data_doc}). 

We define a \emph{mathematician} broadly to include individuals whose doctoral training was completed within departments of mathematics and/or related mathematical sciences. Note that the association of a PhD with the definition of mathematician is a historical one. For example, in \cite{newell1980black} they also focused on mathematicians who earned terminal degrees. Our dataset therefore includes doctoral degrees in theoretical math\footnote{The authors have deliberately chosen use the term 'theoretical mathematics' here in place of the historically used term, 'pure mathematics'. Some mathematicians prefer the term ``theoretical" over ``pure" to describe abstract mathematics because it removes the hierarchical implication that applied work is inherently ``impure" or lesser, instead framing the field as a value-neutral inquiry into rigorous, formal systems and underlying foundational principles.}, applied mathematics, mathematics education, and statistics. Degrees in fields outside of these core areas were considered carefully on a case by case basis in order to best reflect the full range of Black women doctoral degree earners in mathematics.

\subsection{Data Documentation Framework}\label{sec:data_doc}

\begin{step}
    Candidate Mathematician Identification.
\end{step}

 The data for this current project was compiled from analyzing sources including historical scholarship, existing databases and digital archives, institutional and organizational records, and community contributions. Examples include extensive literature searches for journal articles and books documenting participation of Black people in mathematics \cite{MR628027, newell1980black, clark2025hidden}, National Association of Mathematics archives \cite{nam_archive_2026}, archival websites about Black people in science and Mathematics including the MAD \cite{madpages} and MADDER \cite{madder} websites, the Mathematics Genealogy Project \cite{mathgenealogy}, Wikipedia’s List of African American Mathematicians \cite{wiki:african_american_mathematicians}, as well as publicly available lists of women mathematicians such as the list of EDGE Program alumni \cite{edge_program_alumni} and the list of Spelman College alumni in mathematics that have earned PhDs \cite{spelmath_history}. 
 
 In addition to these existing sources, we also identified academic institutions that have been identified as having graduated Black women mathematicians in the past (UC Berkeley, University of Michigan, University of Maryland College Park, Howard University, Emory University, and University of Texas, Austin, Ohio State University) and we reached out to administration in these departments and programs to request graduation records that may be used to identify additional names of Black women that have matriculated from their doctoral programs.   Further, we identified programs (such as the Math Alliance and MSRI-UP) that have been known to include Black women in mathematics at the undergraduate and graduate levels with the goal of supporting them to obtain their doctoral degrees.

\begin{step}
    Candidate Mathematician Verification.
\end{step}

Many of the same sources utilized for identifying candidate mathematicians were also crucial in verifying if the candidate should be in our dataset. 
The crucial task was determining if they met our inclusion criteria for both `Black' and `mathematician'. 
This involved confirming doctoral degree completion, disciplinary classification, potential name variants (many women in our dataset were listed under both maiden and married names), and gender and racial identity across multiple independent sources. 
For a more detailed discussion on some of the processes used in exploring the MADDER database in particular for gender, see Section~\ref{sec:madder}.

\begin{step}
    Metadata collection.
\end{step}

 For each mathematician found, an attempt was made to identify the year each person was born, where they were born, the university that awarded their PhD, the year they earned their degree, their undergraduate institution, and their employment history. We then cross-referenced each data point across multiple sources to verify first, middle, and last names (often last names changed after marriage), verify degree dates, degree granting institutions, dissertation titles, and research areas. For a detailed description of the current metadata see Section~\ref{sec:data_description}.

\begin{step}
Dataset Inclusion and Continuous Curation.
\end{step}

Once a mathematician's records are verified, their data is incorporated into the dataset along with the associated metadata. This data is designed to be a continuously curated resource. There is potential for new Black women to graduate annually with PhDs in mathematics. In addition, new candidates may be identified through additional archival discoveries, institutional records, professional organizations, or community submissions. These future additions will follow the same documentation framework above. Similarly, existing records may be refined as additional information becomes available, including corrections to biographical information, degree records, institutional affiliations, research areas, or other metadata.

\subsection{Data Description}\label{sec:data_description}

The dataset is organized so that each row represents an individual Black woman mathematician and each column records information describing her biographical, educational, and research background. Table~\ref{tab:data_dictionary} summarizes the fields currently maintained in the dataset. We organize the data to include each individual's first name, middle name(s), and last name, together with Boolean variables indicating whether the individual appears in the Mathematicians of the (MADDER) database or the Mathematics Genealogy Project. We additionally record information describing doctoral education, including the year the PhD was awarded, the doctoral degree type (Theoretical Mathematics, Applied Mathematics, Mathematics Education, Statistics, and Other), the PhD granting institution, and the institutional classification (HBCU, PWI, or International), along with up to two doctoral advisors. 

We also collected dissertation titles and classified dissertation subjects by research area using the relevant Mathematical Subject Classification (MSC) codes whenever possible. In addition, the dataset records undergraduate information, including the degree year, degree type, awarding institution, and institutional classification. If an individual earned a master's degree, the same information is collected. When available, geographic and demographic information such as home city, home state, home country, birth and death information was also included.

\begin{table}[htpb]
\centering
\textbf{\normalsize MetaData for the Black Women in Mathematics Dataset} \\[0.75em]
\small
\begin{tabularx}{\textwidth}{p{4cm}p{1.8cm}X}
\toprule
\textbf{Variable} & \textbf{Type} & \textbf{Description / Values} \\
\midrule

\texttt{First\_Name} & String & Given name. \\
\texttt{Middle\_Name\_1} & String & First middle name, when applicable. \\
\texttt{Middle\_Name\_2} & String & Second middle name, when applicable. \\
\texttt{Last\_Name} & String & Family name. \\

\texttt{In\_MAD\_Pages} & Boolean & Indicates whether the individual appears in the MADDER database (Yes/No). \\
\texttt{In\_Math\_Genealogy} & Boolean & Indicates whether the individual appears in the Mathematics Genealogy Project (Yes/No). \\

\texttt{Year\_of\_PhD} & Integer & Year the doctoral degree was awarded. \\
\texttt{Math\_PhD\_Type} & Categorical & Degree classification.
Theoretical Mathematics, Applied Mathematics, Mathematics Education, Statistics, or Other. \\
\texttt{PhD\_Institution} & String & Doctoral degree granting institution. \\
\texttt{PhD\_Institution\_Type} & Categorical &
Institutional classification. Historically Black College or University (HBCU), Predominantly White Institution (PWI), or International. \\

\texttt{PhD\_Advisor\_1} & String & Primary doctoral advisor, when available. \\
\texttt{PhD\_Advisor\_2} & String & Secondary doctoral advisor, when applicable. \\

\texttt{PhD\_Thesis\_Title} & String & Title of the doctoral dissertation. \\
\texttt{PhD\_Research\_Area} & String & Primary research area of the dissertation. \\

\texttt{Bachelor\_Year} & Integer & Year bachelor's degree was awarded. \\
\texttt{Bachelor\_Type} & Categorical & Degree classification. Theoretical Mathematics, Applied Mathematics, Mathematics Education, Statistics, or Other.  \\
\texttt{Bachelor\_Institution} & String & Bachelor's degree-granting institution. \\
\texttt{Bachelor\_Institution\_Type} & Categorical & Institutional classification. Historically Black College or University (HBCU), Predominantly White Institution (PWI), or International. \\

\texttt{Master\_Year} & Integer & Year master's degree was awarded, when applicable. \\
\texttt{Master\_Type} & Categorical & Degree classification. Theoretical Mathematics, Applied Mathematics, Mathematics Education, Statistics, or Other. \\
\texttt{Master\_Institution} & String & Master's degree-granting institution. \\
\texttt{Master\_Institution\_Type} & Categorical & Institutional classification. Historically Black College or University (HBCU), Predominantly White Institution (PWI), or International, when applicable. \\

\texttt{Home\_City} & String & Home city, when available. \\
\texttt{Home\_State} & String & Home state or province, when available. \\
\texttt{Home\_Country} & String & Home country. \\

\texttt{Birth} & Date & Date or year of birth, when available. \\
\texttt{Death} & Date & Date or year of death, when applicable. \\

\bottomrule
\end{tabularx}

\caption{Data dictionary for the Black Women in Mathematics dataset. Each row of the dataset corresponds to one individual.}

\label{tab:data_dictionary}
\end{table}

\subsection{MADDER}\label{sec:madder}

In order to obtain how many of the 502 Black PhDs in the mathematical sciences listed in the MADDER database should be included in our set of Black women mathematicians, we conducted a systematic search of publicly available information. 
We first examined each individual's MAD page for an available photograph and biographical information, as biographies may explicitly identify an individual's pronouns. When neither a photograph nor a biography was available on the MAD page, we conducted a broader Google search for publicly accessible information. 
Sources consulted included faculty and institutional websites, professional webpages including but not limited to LinkedIn profiles, newspaper articles, obituaries, and RateMyProfessors profiles that explicitly identified an individual's pronouns. 

\begin{figure}[htpb]
\centering
\textbf{\normalsize Black Math PhDs in MADDER Characterized by Gender by Decade} \\[0.75em]
 \includegraphics[width=\textwidth]{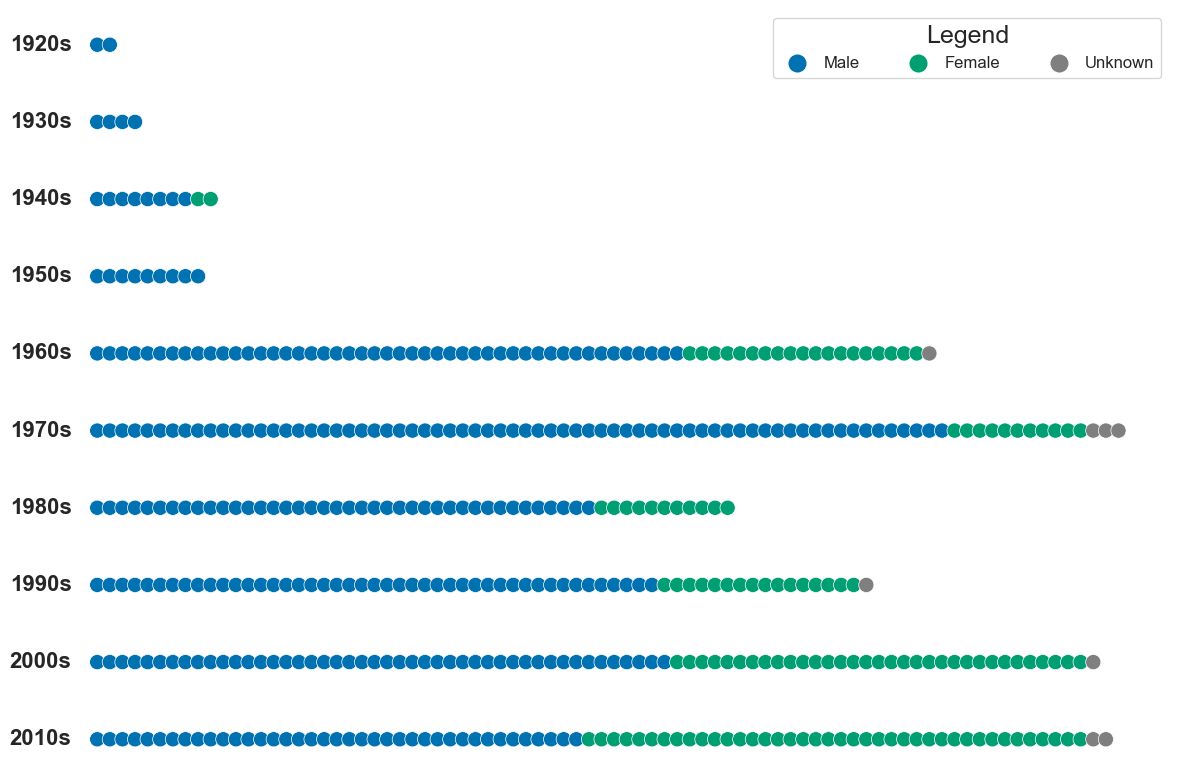}
    \caption{The MADDER dataset organized by gender and by decade.}
    \label{fig:madder_gender}
\end{figure}

For the vast majority of individuals, pronouns were identified through photographs and/or explicit pronoun usage in publicly available biographical or professional sources. 
For 4 individuals, pronouns were inferred from their names when no other publicly available information was available. 
We were unable to determine the pronouns of 14 individuals. For the purposes of this study, individuals using she/her pronouns were categorized as female, those using he/him pronouns were categorized as male, one individual used she/they pronouns and was categorized as other. \footnote{Gender was inferred from publicly available pronoun use, which may not accurately reflect individuals' gender identities.} Among those for whom pronouns could be identified, 145 individuals used she/her pronouns, 331 used he/him pronouns, and one used she/they pronouns. 
During the process, we also identified eight duplicate records in the database. Figure~\ref{fig:madder_gender} shows the result of this work where Black Math PhDs in MADDER are characterized by gender and by decade.

\subsection{Limitations}

Although every effort was made to collect complete information for each individual, not every variable is available for every record. In some instances the lack of availability of public and institutional records limited our ability to collect complete metadata for each individual. 
In these cases, rather than exclude the individual from the dataset entirely, they were kept if we verified that they met our baseline inclusion criteria of having a PhD in the mathematical sciences and were from the African Diaspora. The dataset is designed to be updated as additional individuals are identified and existing records are refined. 

The methodology described in this section provides a transparent and reproducible framework for documenting Black women mathematicians in the United States. While the dataset will continue to evolve as additional information becomes available, its current scope is sufficient to support an initial quantitative examination of educational pathways, advisors, institutional representation, and research contributions. The following section presents exploratory analyses that illustrate these emerging patterns.


\section{Black Women Mathematicians in the United States: Findings}\label{sec:results}

 Using the curated dataset described in Section~\ref{sec:methods}, this section presents the first quantitative portrait of Black women mathematicians in the United States through a series of tables and visualizations. These descriptive analyses examine patterns in doctoral production, undergraduate preparation, advisor networks, institutional affiliations, and research areas. The main results show that doctoral production is concentrated at a relatively small number of institutions and that Historically Black Colleges and Universities (HBCUs) play a central role in undergraduate preparation. Furthermore, doctoral guidance is broadly distributed, but includes several particularly influential advisors. Finally, the mathematical contributions of Black women span multiple research areas with substantial representation in applied mathematics and mathematics education. These findings represent the current state of the curated dataset and will continue to evolve as additional records are identified, verified, and incorporated.

\subsection{Doctoral Institutions}

The institutions represented in our dataset illustrate both the breadth and the concentration of doctoral training for Black women mathematicians. Of the more than 500 doctoral programs currently identified as offering degrees in statistics, theoretical mathematics, or applied mathematics in the United States and Canada, only 63 are represented in our dataset\footnote{Note that 3 women in our dataset received their PhDs from institutions outside of North America.}. In other words, based on currently available records, approximately 13\% of institutions offering doctoral degrees in the mathematical sciences are known to have awarded a PhD in mathematics to a Black woman. 

\begin{figure}[htpb]
\centering
\textbf{\normalsize Top 10 PhD Granting Institutions of Black Women in Mathematics in our Dataset and their Proportions to all Represented Institutions} \\[0.5em]
\includegraphics[width=0.95\textwidth]{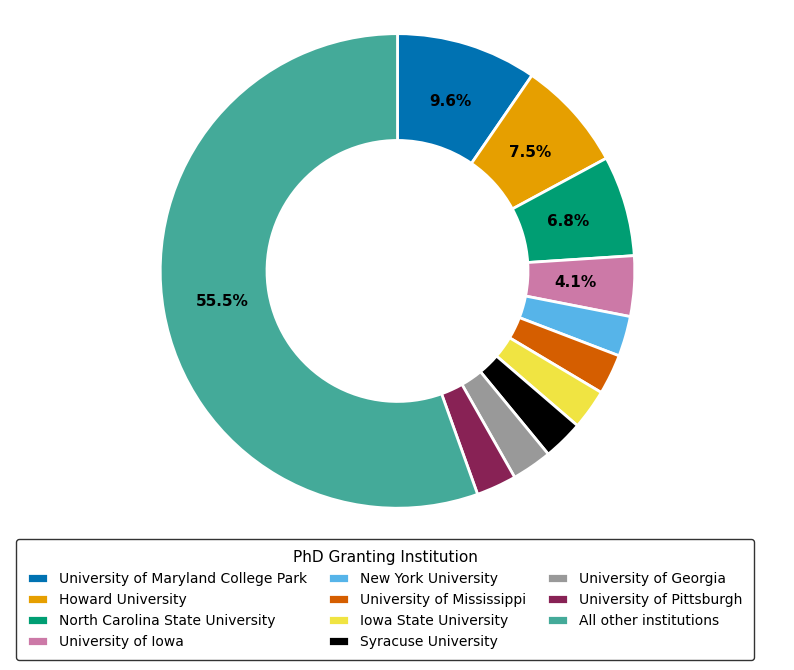}
    \caption{Percentage of PhDs awarded to Black women by University. The top 10 PhD institutions with respect number of Black women mathematicians graduated are indicated in the legend. They account for 44.5\% of all institutions in our data, almost half.}
\label{fig:PhD_Institutions_circle}
\end{figure}

Although doctoral training has occurred across dozens of universities, a small number of institutions have played especially prominent roles. The University of Maryland, College Park (14 graduates) and Howard University (11 graduates\footnote{At the time of this article, we are aware of additional Black women PhDs who have earned degrees from Howard. Some of these names have already been identified, and the next step is verification and the collection of metadata. Total degree production from Howard, and other HBCUs with PhD programs, Delaware State and Morgan State, can be searched for through graduate commencements and institutional data requests.}) have awarded the largest numbers of doctoral degrees to Black women mathematicians, followed by North Carolina State University (10) and the University of Iowa (6). The remaining institutions each produced four or fewer graduates, underscoring both the rarity of these degrees and the broad distribution of doctoral mentorship across the mathematical sciences. These concentrations suggest that doctoral production among Black women mathematicians has not occurred uniformly across institutions. Of the institutions that have awarded a math PhD to Black women, the majority of them have only done so once. Institutions that awarded doctoral degrees to three or more Black women mathematicians are listed in Table \ref{tab:institution_summary_table} and visualized in Figure \ref{fig:PhD_Institutions_circle}. 

Another story emerges when undergraduate and doctoral institutions are compared. Currently, institutions that awarded Black women the most Bachelor degrees are HBCUs, whereas institutions that have awarded Black women the most math PhDs are Predominantly White Institutions (PWIs), with the exception of Howard University. This makes some sense historically as Howard was the first HBCU to offer a PhD in mathematics in 1976, a full 33 years after the first Black woman earned a PhD in mathematics in 1943. Nonetheless, the current lack of comparable concentration among predominantly white institutions at the undergraduate level is an intriguing question for future research.

\begin{table}[!ht]
\centering
\textbf{\normalsize Institutions who have Awarded Degrees to \\ Three or More Black Women Mathematicians} \\[0.75em]
\begin{tabular}{lrr}
\toprule
\textbf{Institution} & 
\textbf{Bachelor Degree} & 
\textbf{Doctoral Degree} \\
\midrule
Spelman College & 27 & -- \\
Howard University & 5 & 11 \\
Fisk University & 4 & -- \\
Hampton University & 4 & -- \\
Morgan State University & 4 & -- \\
Xavier University of Louisiana & 4 & -- \\
Alcorn State University & 3 & -- \\
Clark Atlanta University & 3 & -- \\
Dillard University & 3 & -- \\
North Carolina A\&T State University & 3 & -- \\
South Carolina State University & 3 & -- \\
University of Maryland, College Park & -- & 14 \\
North Carolina State University & -- & 10 \\
University of Iowa & -- & 6 \\
Iowa State University & -- & 4 \\
New York University & -- & 4 \\
Syracuse University & -- & 4 \\
University of Georgia & -- & 4 \\
University of Mississippi & -- & 4 \\
University of Pittsburgh & -- & 4 \\
City University of New York & -- & 3 \\
Emory University & -- & 3 \\
Georgia State University & -- & 3 \\
Rice University & -- & 3 \\
University of Alabama & -- & 3 \\
University of California, Berkeley & -- & 3 \\
University of Michigan & -- & 3 \\
University of Nebraska--Lincoln & -- & 3 \\
The Ohio State University & -- & 3 \\
\bottomrule
\end{tabular}

\caption{Institutions who awarded degrees to three or more Black women mathematicians in the dataset. Undergraduate and doctoral institutions are shown separately. Institutions with unavailable educational information or whose records are still under review are omitted.}
\label{tab:institution_summary_table}
\end{table}

Some scholarship already exists exploring institutional patterns of Black women in mathematics. In 2022 the article ``An Unlikely Cohort: Disrupting the Stereotypes about Mathematicians and their Careers” was published in the \textit{Notices of the AMS} \cite{MR4378418}. 
There, Taryn Butler Lewis, Tasha R. Inniss, Monica Jackson, and Calandra Tate Moore told the story of the 20 Black Women who have either a Master’s or Doctoral degree or both in Mathematics from the University of Maryland College Park. 
This was followed by the article ``Spelman College: a model of success for producing black women in mathematics” in 2023 in the Notices of the American Mathematical Society \cite{MR4537174}, highlighting the success that the Mathematics Department at Spelman College has had in training Black Women for careers as mathematicians.'' 

\begin{figure}[htpb]
    \centering
    \textbf{\normalsize Migration Patterns to Doctoral Institutions from Black Women PhDs who Earned Undergraduate Degrees from Spelman College} \\[0.75em]
    \includegraphics[width=\textwidth]{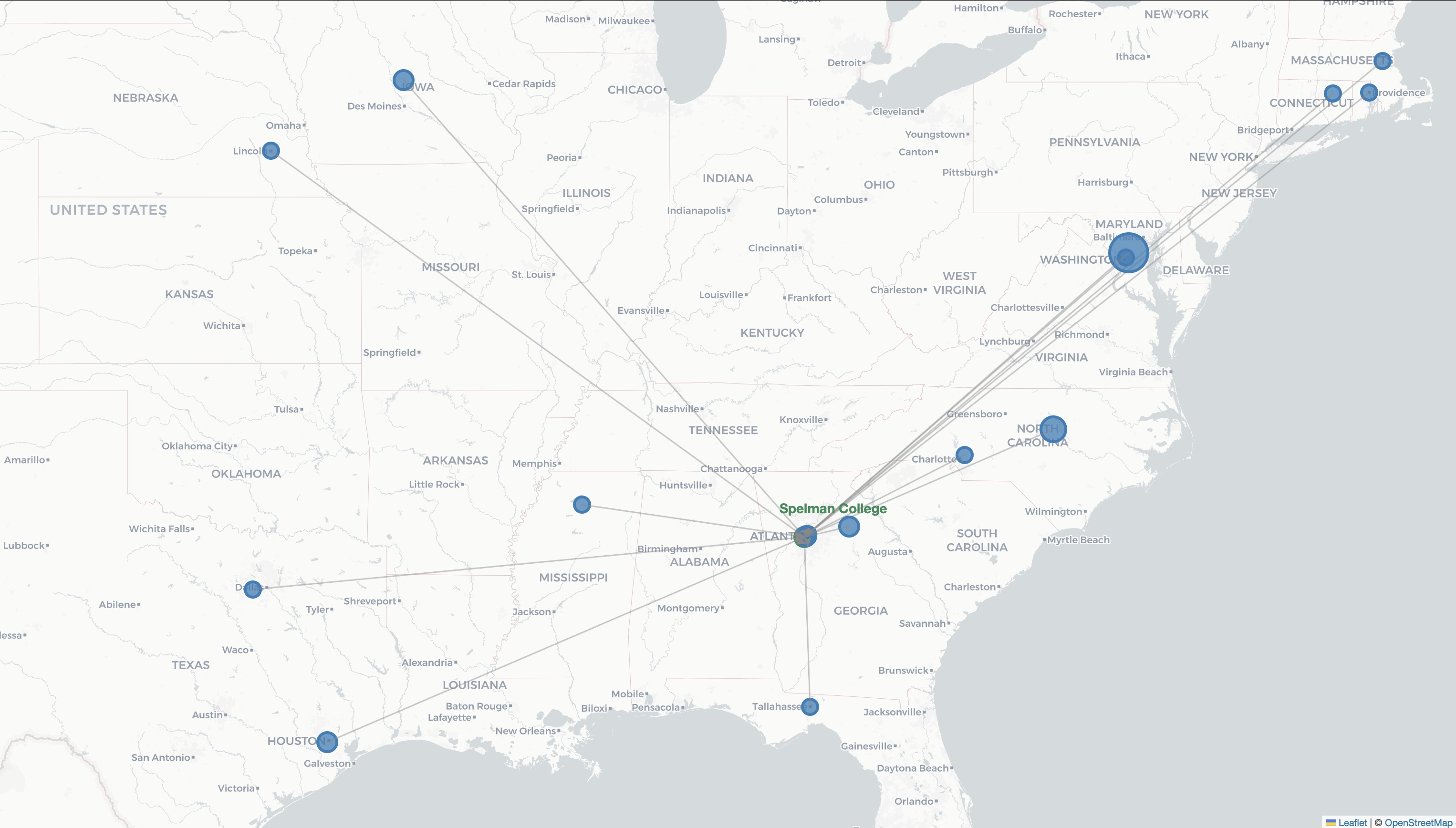}
    \caption{The lines emanating from Atlanta, Georgia show pathways from one undergraduate institution to doctoral institutions for the 27 Black Women PhDs who earned undergraduate degrees in math from Spelman in our current dataset. These include the University of Maryland College Park (6), North Carolina State University (6), the University of Georgia (2), Rice University(2), Iowa State University (2), Georgia State University (2), the University of Nebraska Lincoln (1), Massachusetts Institute of Technology (1), Howard University (1), University of North Carolina Charlotte (1), University of Mississippi (1), Brown University (1), Florida State University (1), University of Connecticut (1), University of Texas Arlington (1), and Emory University (1).}
    \label{fig:spelman_migration}
\end{figure}

Figure~\ref{fig:spelman_migration} highlights migration patterns of Black women in mathematics who have earned Bachelor degrees from Spelman college. The lines emanating from Spelman college indicate key geographical areas where doctoral degrees where earned. The large circle near Washington D.C., for example, indicates the many Black women from Spelman who earned their doctorates from the University of Maryland, College Park.

\subsubsection{Research Areas}

In addition to documenting who has earned PhDs in mathematics, our dataset provides an opportunity to examine the intellectual breadth of Black women's contributions to the discipline. We classified each dissertation according to its primary and secondary Mathematics Subject Classification (MSC) codes. This was done manually through review of dissertation titles and abstracts and then verified with the assistance of a large language model which was used as a tool for suggestion, organization, and verification of classifications. This enabled us to identify broad patterns in the mathematical areas in which Black women have conducted doctoral research. 

\begin{figure}[htpb]
    \centering
    \textbf{\normalsize Cumulative Growth of Top Five Research Areas \\ for Black Women Mathematicians by Decade} \\[0.75em]
 \includegraphics[width=\textwidth]{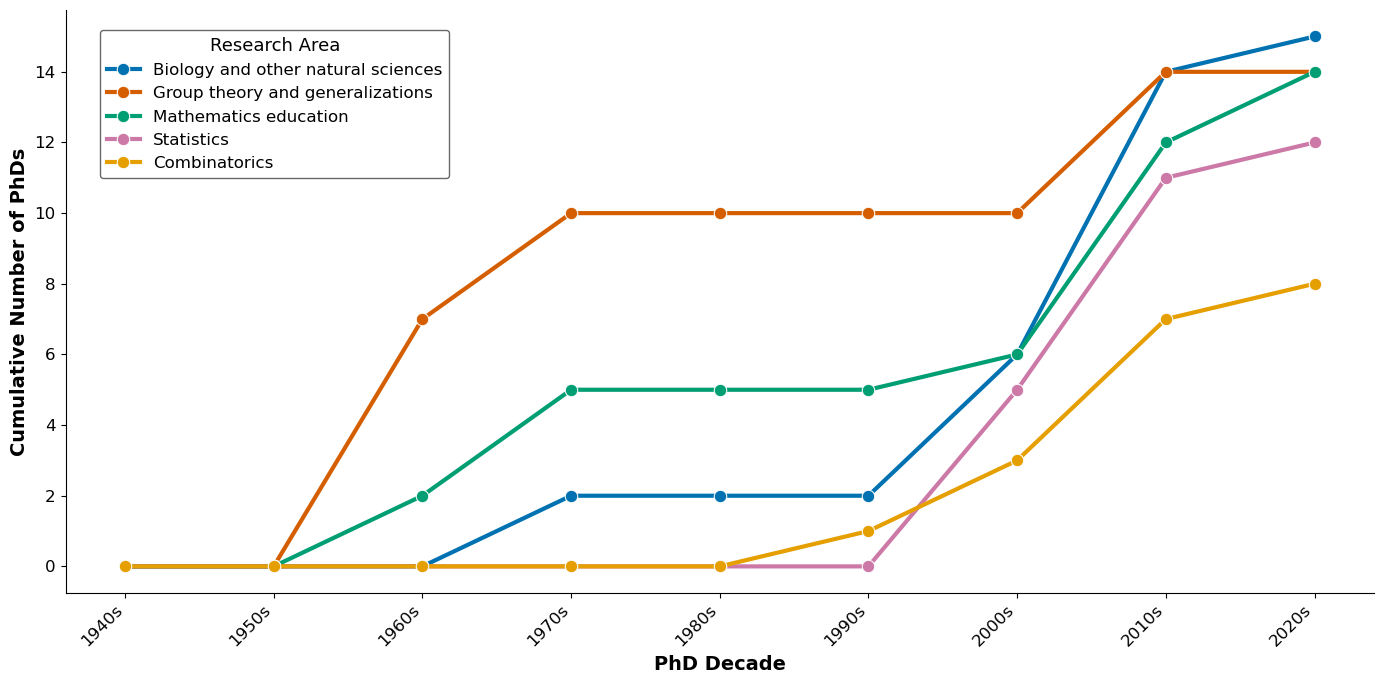}
    \caption{The top five research areas according to primary MSC classification are highlighted in the legend above: biology and other natural science, group theory and generalizations, mathematics education, statistics, and combinatorics. Cumulative trends for each of these areas are shown by decade, starting in 1940.}
    \label{fig:cumulative_research}
\end{figure}

Our first exploration sought to understand primary research areas that were common among Black women mathematicians. The cumulative growth of the five most common primary dissertation
research areas among Black women mathematicians since 1940 is shown in Figure~\ref{fig:cumulative_research}.
When considering both primary and secondary MSC classifications, the most common pairing of MSC classifications was Biology and Natural Sciences as the primary classification and Ordinary Differential Equations as the secondary classification, showcasing a strong interdisciplinary scholarship among these women. Mathematics Education emerged as the second most common primary research area. Notably, these two areas, one centered on the application of mathematics to problems in the natural sciences and the other on mathematics teaching and learning, have both been identified in the literature as fields in which underrepresented students in STEM have historically found important opportunities for participation and persistence \cite{goodwin2024represented, echegoyen2025factors, ashley2017building, stephens2017undergraduate}.

\subsubsection{PhD Advisors}
This section recognizes and documents the faculty members who advised Black Women through their doctoral journeys. Earning a PhD is never a solitary endeavor, and thoughtful mentorship often plays an essential role in helping students develop as scholars, navigate the challenges of graduate education, and launch successful research careers. Our extensive research identified 135 unique PhD advisors among the Black women mathematicians in our dataset. Given that the count of our data set at the time of this article is 147, this implies that the majority of these faculty members have only supervised one Black woman PhD advisee, with a very small subset of faculty advising two. This distribution underscores that the responsibility for mentoring Black women mathematicians has historically been spread across many individuals rather than concentrated among a small number of faculty members. A list of PhD advisors who supervised more than one Black woman mathematician is outlined in Table \ref{tab:multiple_advisors}.

\begin{table}[htpb]
\centering
\textbf{\normalsize Advisors who have Supervised Two or More Black Women PhDs} \\[0.75em]
\begin{tabular}{llr}
\toprule
\textbf{PhD Advisor} & \textbf{Institution} &  \textbf{Number of Advisees} \\
\midrule
Neil Hindman &Howard University & 9 \\
James Dolan Reid & Syracuse University & 4 \\
Alan Hopenwasser & University of Alabama & 2 \\
Clinton G. McCrory & University of Georgia & 2 \\
Doron Levy & University of Maryland College Park & 2 \\
Edward DeRoche & Marquette University & 2 \\
Frank Charles Hoppensteadt & New York University & 2 \\
Gilbert Baumslag & City University of New York & 2 \\
Mark William Brittenham & University of Nebraska Lincoln & 2 \\
Mette Olufsen & North Carolina State University & 2 \\
Norman Levine & University of Pittsburgh & 2 \\
Richard Scott Pierce & University of Washington & 2 \\
Talmage James Reid & University of Mississippi & 2 \\
W. E. Deskins & University of Pittsburgh & 2 \\
William H. Kazez & University of Georgia & 2 \\
\bottomrule
\end{tabular}
\caption{PhD advisors who supervised more than one Black woman mathematician in the dataset.}
\label{tab:multiple_advisors}
\end{table}

One individual stands out for the sustained mentorship they provided. Dr. Neil Hindman, a professor of mathematics at Howard University, advised nine Black women mathematicians. In 2006 the MAD Pages noted, ``that Howard University Professor Neil Hindman has been Ph.D. advisor to more African Americans than anyone else in the world.'' While this project cannot attest to this claim for all Black mathematicians, Dr. Hindman has supervised more Black women in the mathematical sciences to completion than any other advisor identified in our dataset. Find a detailed account on the work of Dr. Hindman's students in this article, ``The Research of Thirteen Students at Howard University", written by Dr. Hindman himself \cite{hindman2009research}. 
We recognize Dr. Hindman not only for the number of dissertations he supervised, but for the meaningful role he played in supporting multiple Black women through one of the most consequential stages of their mathematical careers.

\subsection{Undergraduate Institutions and the Role of HBCUs}

Thus far we have discussed institutions that awarded doctoral degrees to Black women in mathematics, but it is also informative to examine where these mathematicians began their academic journeys. Undergraduate institutions play a formative role in preparing students for graduate study by providing early research experiences, mentoring, and encouragement to pursue advanced degrees. We have not yet prioritized similar analyses of master's degree granting institutions due to both incomplete data as well as the fact that many mathematicians proceed directly from a bachelor's degree to a doctoral program without earning a terminal master's degree. Table~\ref{tab:institution_summary_table} lists the undergraduate institutions attended by three or more Black women represented in our dataset. Note Table \ref{tab:institution_summary_table} excludes 7 women for whom we could not find any information on their bachelor's institution.

The undergraduate institutions represented in our dataset reveal a striking pattern. Every institution in our dataset that awarded bachelor's degrees to three or more Black women who later earned doctorates in mathematics is a Historically Black College or University (HBCU). This finding underscores the central role that HBCUs have played in cultivating generations of Black women mathematicians and preparing them for doctoral study. Among these institutions, Spelman College stands out as the clear leader, with 27 alumnae in our dataset who went on to earn PhDs in mathematics (Table ~\ref{tab:institution_summary_table}). Figure~\ref{fig:spelman_migration} illustrates the doctoral pathways
of Black women mathematicians who earned their bachelor's degrees
from Spelman College.
Howard University follows with five graduates, while Fisk University, Hampton University, Morgan State University, and Xavier University of Louisiana each contributed four. 

\begin{figure}[h]
    \centering
    \textbf{Black Women PhDs in Mathematics by Decade and HBCU Undergraduate Classification}
    \includegraphics[width=\textwidth]{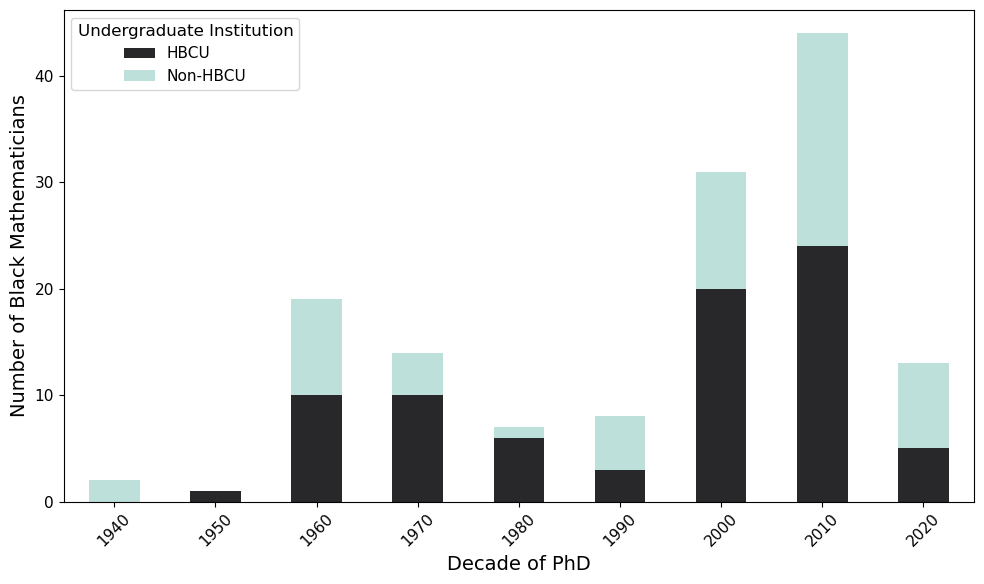}
    \caption{The number of Black women earning a PhD in mathematics from 1940 to 2020 is show by decade. Each stacked bar in the plot indicates how many of the total PhDs that decade earned undergraduate degrees at an HBCU.}
    \label{fig:bw_phd_decade}
\end{figure}

In addition, Figure~\ref{fig:bw_phd_decade} illustrates these trends by decade. For example, according to our current dataset, in the 1950s and in the 1980s, almost all Black women earning a PhD in the mathematical sciences came from an undergraduate HBCU. Even in the 2010s, it appears that half of Black women PhD receipients in mathematicians earned undergraduate degrees from an HBCU. The prominence of HBCUs is consistent with decades of research documenting their outsized contribution to the education of Black scientists, engineers, and mathematicians despite representing only a small fraction of colleges and universities in the United States. \\

\section{Conclusion}\label{sec:conclusion}

The current dataset establishes a foundation for a broad program of research on Black mathematicians. While the primary focus of this project is the documentation of women mathematicians, the creation of this dataset instituted standardized archival processes that can be expanded to include all Black mathematicians and other groups of interest. In our conclusion, we discuss avenues for future work with an emphasis on how this effort created opportunities for both quantitative and historical investigations on this topic that have previously been difficult or impossible due to the absence of systematic data. We also note the connection between this work and the growing mathematical research community in quantitative justice (Section~\ref{sec:mm_qj}).

\subsection{Future Work}\label{sec:futurework}

In future work, our emphasis will shift from Black women specifically, to any person from the African diaspora that meets the aforementioned criteria of \emph{mathematician}, regardless of gender. Below we describe a broad range of future investigations broken down thematically.

\begin{theme}\label{theme:doc}
Historical documentation.
\end{theme}

We will solicit correspondence to make it simple for individuals to (self) report Black mathematicians. 
We have already done this through newsletter announcements \cite{awm2026enumerating, nam2026phdsdiaspora} and a shareable google form. 
We encourage anyone who has information about Black mathematicians to reach out to our team\footnote{We provide a link to a form for (self) reporting Black mathematicians: \url{https://tinyurl.com/psv6cm2e}.}. 
Furthermore, as we continue documenting Black mathematicians, we will continuously improve our methodology. 
This includes using archival methods to collect university commencement records, graduate school archives, dissertation databases, mathematics department records, professional association archives, obituaries, and biographical collections.

\begin{theme}\label{theme:collab}
Collaboration with existing documentation efforts.
\end{theme}

This is closely related to Research Theme~\ref{theme:doc} but distinct enough to warrant its own discussion. Our dataset overlaps with the broader MADDER project. We will collaborate to update this database. In addition, during our documentation process, it was noted that many Black women were missing from the Math Genealogy\footnote{The Mathematics Genealogy Project catalogs over 330,000 individuals worldwide, but it is self-reported and self-corrected; its coverage of underrepresented scholars has been uneven since its founding \citep{jackson2007, mulcahy2017}.} project. This provides another avenue for potential collaboration.  

We may also collaborate on initiatives for documenting underrepresented mathematicians in general. There are several digital initiatives to document and celebrate mathematicians from historically underrepresented groups through public facing online archives. \emph{Lathisms} \cite{lathisms}, \emph{Mathematically Gifted and Black} \cite{GiftedandBlack}, and \emph{Indigenous Mathematicians} \cite{indigenousmathematicians}, are all websites showcasing contributions from Hispanic and Latinx, Black, and Indigenous Mathematicians respectively\footnote{\emph{Lathisms} was founded in 2016 by Alexander D\'iaz-L\'opez, Pamela Harris, Alicia Prieto Langarica, and Gabriel Sosa \cite{lathisms}. Each year during Hispanic Heritage Month, daily biographical profiles are posted with  information that includes relevant research, teaching, and service contributions of their honorees. Later in 2016, this website along with the Hidden Figures major motion picture, inspired another website, \emph{Mathematically Gifted and Black}, which featured Black Mathematicians and daily profiles released during Black History Month. This was founded by Erica Graham, Raegan Higgins, Shelby Wilson, and Candice Price \cite{GiftedandBlack}. In 2019, a website spearheaded by a `Pacific Islanders in Mathematics' session at the Society for the Advancement of Native Americans and Chicanos in the Sciences (SACNAS) called \emph{Indigenous Mathematicians} was created. This site aims to archive all known information about Indigenous Mathematicians across the world, and also features profiles and podcast of select honorees. \cite{indigenousmathematicians}.}. Each features biographical profiles of their honorees. Holistic documentation of these communities is also essential to understanding historical patterns of participation in mathematics.   

\begin{theme}
Doctoral production and institutional pathways.
\end{theme}

The dataset makes it possible to examine long term trends in doctoral production by Black mathematicians. As the dataset grows, it will become possible to study how the number of Black mathematicians earning doctoral degrees in mathematics and related fields has changed over time, identify periods of growth or stagnation, and explore the institutional and societal factors associated with those changes. Such analyses can contribute to a more complete historical understanding of factors that impact participation in the mathematical sciences.

In addition, we can identify universities and departments that have played significant roles in educating Black mathematicians. Examining patterns of doctoral production across institutions may reveal important insights into recruitment, retention, mentoring, and educational pathways. Similarly, geographic analyses can help illuminate patterns of mobility, regional concentrations of mathematical talent, and the relationships between undergraduate and doctoral institutions.

\begin{theme}
Academic genealogy and mentorship.
\end{theme}

By integrating this dataset with resources such as the Mathematics Genealogy Project (see Research Theme~\ref{theme:collab}), researchers may be able to trace advising relationships, identify influential mentoring communities, and better understand how knowledge, opportunity, and professional support are transmitted across generations of mathematicians. Such work would provide a richer picture of the social and intellectual structures that shape mathematical careers.

\begin{theme}
Career trajectories and professional outcomes.
\end{theme}

As additional information is collected, the dataset may also support investigations of career trajectories, research specializations, leadership roles, and professional outcomes. These analyses could help document the diverse contributions Black women have made across academia, industry, government, education, and other sectors. Understanding these pathways has the potential to inform future mentoring initiatives, professional development programs, and efforts to broaden participation in the mathematical sciences.\\

\subsection{Metamath and Quantitative Justice}\label{sec:mm_qj}

The research themes described in the previous section illustrate how historical documentation can serve as the foundation for rigorous mathematical research. The 2025 AIM workshop on \emph{MetaMath} is a proof of concept for this perspective, as it brought together researchers interested in using mathematics, statistics, and data science to study the mathematical sciences community itself. Viewed through this lens, the dataset presented in this article can be seen as research infrastructure that enables quantitative investigation of doctoral production, institutional pathways, academic genealogy, mentorship networks, geographic mobility, and professional trajectories. These questions can be studied using a broad range of mathematical and computational techniques, including statistical analysis, network science, graph theory, machine learning, data visualization, and other methods from the mathematical sciences.

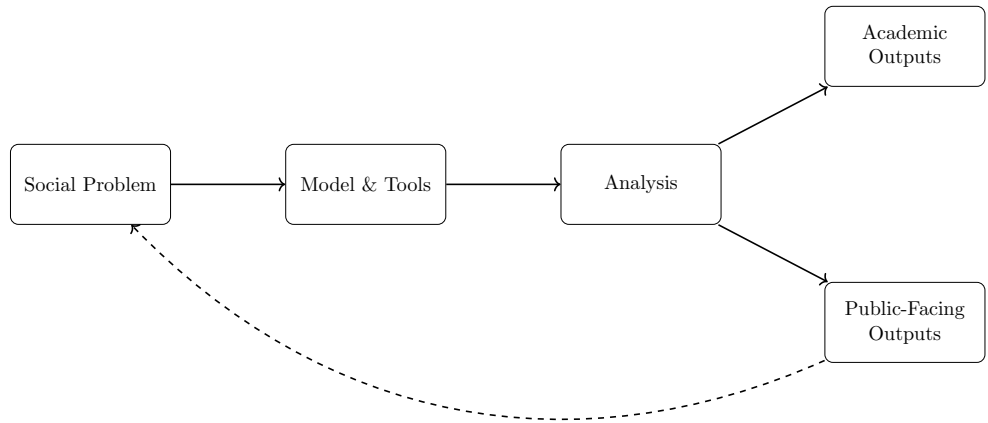
\begin{figure}[h]
\centering
\textbf{\normalsize Framework for Quantitative Justice Research in Mathematics} \\ [0.75em]
\resizebox{\textwidth}{!}{
\begin{tikzpicture}[scale = 1.2,
    node distance=2cm,
    every node/.style={
        draw, 
        rectangle, 
        rounded corners, 
        align=center, 
        minimum width=2.8cm, 
        minimum height=1.4cm,
        font = \normalsize
    },
    arrow/.style={->, thick}
]

\node (input) {Social Problem};
\node (model) [right=of input] {Model \& Tools};
\node (analysis) [right=of model] {Analysis};

\node (academic) [above right=1cm and 1.8cm of analysis] {Academic \\ Outputs};
\node (community) [below right=1cm and 1.8cm of analysis] {Public-Facing \\ Outputs};

\draw[arrow] (input) -- (model);
\draw[arrow] (model) -- (analysis);
\draw[arrow] (analysis) -- (academic);
\draw[arrow] (analysis) -- (community);
\draw[->, dashed, thick, bend left=35] (community) to (input);

\end{tikzpicture}
}
\caption{A conceptual framework for quantitative justice research. Social questions motivate the selection of mathematical models and analytical tools, producing both traditional scholarly outputs and broader community-oriented impact. The dashed feedback arrow represents the iterative relationship between societal outcomes and future mathematical inquiry.}  
\label{fig:qj_framework}
\end{figure}

More broadly, this work illustrates one way in which MetaMath contributes to the emerging field of quantitative justice. Quantitative justice seeks to address societal issues related to inequity with rigorous theoretical and applied mathematical tools. 
In particular, its framework, as illustrated in Figure~\ref{fig:qj_framework}, begins with societal questions and then develops rigorous mathematical methods to investigate them. 
The societal issue motivating this project concerns the historical participation and representation of Black mathematicians. 
Our project highlights that the counting and documentation described in this article can be used as a tool to facilitate models to analyze the mathematical community. 

In turn, our growing dataset of 147 Black women mathematicians provides a foundation for future mathematical research while simultaneously preserving the history of Black mathematicians and making that history more accessible to researchers, educators, professional societies, and the broader public. 
Thus our work serves not only as a foundation for MetaMath research on the mathematical sciences community, but can be understood in broader contexts as a novel approach to quantitative justice research that models how archival scholarship, mathematics, statistics, and data science can work together to address broader questions of equity, participation, and representation.

\section{Declaration of Competing Interests}
On behalf of all authors, the corresponding author states that there is no conflict of interest.


\newpage

\appendix

\renewcommand{\thesection}{Appendix~\Alph{section}}

\section{Bachelor's Degree Institutions}

This is a list of undergraduate institutions that have granted Bachelors degrees to Black women who later earned PhD in mathematical sciences, ordered by the number of degree recipients (in descending order) and alphabetically in the case of ties. The data was compiled from MAD Pages, professional websites, newspaper articles, doctoral dissertations, and publicly available sources.

\begin{longtable}{p{0.78\textwidth}r}
\caption{Undergraduate institutions attended by Black women mathematicians represented in the dataset.}
\label{tab:full_bachelors_institutions}\\
\toprule
\textbf{Undergraduate Institution} & \textbf{Number of Graduates} \\ \midrule
\endfirsthead

\multicolumn{2}{c}{\tablename\ \thetable\ -- continued from previous page}\\
\toprule
\textbf{Undergraduate Institution} & \textbf{Number of Graduates} \\ \midrule
\endhead

\midrule
\multicolumn{2}{r}{Continued on next page}\\
\endfoot

\bottomrule
\endlastfoot

Spelman College & 27 \\
Howard University & 5 \\
Fisk University & 4 \\
Hampton University & 4 \\
Morgan State University & 4 \\
Xavier University of Louisiana & 4 \\
Alcorn State University & 3 \\
Bryn Mawr College & 3 \\
Clark Atlanta University & 3 \\
Dillard University & 3 \\
North Carolina A\&T State University & 3 \\
South Carolina State University & 3 \\
Alabama A\&M University & 2 \\
Florida A\&M University & 2 \\
Grambling State University & 2 \\
Smith College & 2 \\
Southern University & 2 \\
Texas Southern University & 2 \\
University of Kansas & 2 \\
University of Kentucky & 2 \\
University of Mississippi & 2 \\
Antioch College & 1 \\
Baylor University & 1 \\
Brandeis University & 1 \\
Bucknell University & 1 \\
California State University, Chico & 1 \\
California State University, Long Beach & 1 \\
City University of New York & 1 \\
City University of New York - College of Staten Island & 1 \\
College of Bahamas, Bahamas & 1 \\
Columbia College of South Carolina & 1 \\
Florida Atlantic University & 1 \\
Harvard University & 1 \\
Hofstra University & 1 \\
Hood College & 1 \\
Hunter College & 1 \\
Indiana University Bloomington & 1 \\
Johnson C. Smith College & 1 \\
LeMoyne College & 1 \\
Marianao Institute & 1 \\
Mercer University & 1 \\
Miami University & 1 \\
Mount Holyoke College & 1 \\
Northwestern University & 1 \\
Oakwood University & 1 \\
Paine College & 1 \\
Pomona College & 1 \\
Rhodes College & 1 \\
Rose-Hulman Institute of Technology & 1 \\
Sam Houston State University & 1 \\
Seton Hall University & 1 \\
Shaw University & 1 \\
Southern University at New Orleans & 1 \\
Stephen F. Austin University & 1 \\
Temple University & 1 \\
The College of New Jersey & 1 \\
The Ohio State University & 1 \\
University of Alabama & 1 \\
University of Arkansas at Pine Bluff & 1 \\
University of California, Berkeley & 1 \\
University of Chicago & 1 \\
University of Florida & 1 \\
University of Georgia & 1 \\
University of Ghana & 1 \\
University of Hawaii at Hilo & 1 \\
University of Ibadan, Nigeria & 1 \\
University of Ilorin, Nigeria & 1 \\
University of Maryland, College Park & 1 \\
University of Michigan & 1 \\
University of Pittsburgh & 1 \\
University of Rochester & 1 \\
University of South Carolina & 1 \\
University of Texas at Arlington & 1 \\
Vassar College & 1 \\
Williams College & 1 \\
Winston-Salem State University & 1 \\
\end{longtable}

\noindent\textit{Note.} Undergraduate institution information was unavailable for 7 women in the dataset.

\vspace{0.5in}

\section{Doctoral Degree Institutions}

This is a list of doctoral institutions that awarded Ph.D.s in mathematical sciences to Black women, ordered by the number of degree recipients (in descending order) and alphabetically in the case of ties. The data was compiled from The Math Geneology Project, MAD Pages, professional websites, newspaper articles, doctoral dissertations, and publicly available sources.

\begin{longtable}{p{0.78\textwidth}r}
\caption{PhD-granting institutions represented in the dataset.}
\label{tab:full_phd_institutions}\\
\toprule
\textbf{PhD-Granting Institution} & \textbf{Number of PhDs} \\
\midrule
\endfirsthead

\multicolumn{2}{c}{\tablename\ \thetable\ -- continued from previous page}\\
\toprule
\textbf{PhD-Granting Institution} & \textbf{Number of PhDs} \\
\midrule
\endhead

\midrule
\multicolumn{2}{r}{Continued on next page}\\
\endfoot

\bottomrule
\endlastfoot

University of Maryland, College Park & 14 \\
Howard University & 11 \\
North Carolina State University & 10 \\
University of Iowa & 6 \\
Iowa State University & 4 \\
New York University & 4 \\
Syracuse University & 4 \\
University of Georgia & 4 \\
University of Mississippi & 4 \\
University of Pittsburgh & 4 \\
City University of New York & 3 \\
Emory University & 3 \\
Georgia State University & 3 \\
Rice University & 3 \\
University of Alabama & 3 \\
University of California, Berkeley & 3 \\
University of Michigan & 3 \\
University of Nebraska--Lincoln & 3 \\
University of Texas at Arlington & 3 \\
Marquette University & 2 \\
Michigan State University & 2 \\
Temple University & 2 \\
The Ohio State University & 2 \\
University of Alabama at Birmingham & 2 \\
University of Massachusetts Amherst & 2 \\
University of Utah & 2 \\
University of Washington & 2 \\
Virginia Commonwealth University & 2 \\
Auburn University & 1 \\
Bielefeld University, Germany & 1 \\
Binghamton University, SUNY & 1 \\
Brown University & 1 \\
The Catholic University of America & 1 \\
Cornell University & 1 \\
Florida Atlantic University & 1 \\
Florida State University & 1 \\
The George Washington University & 1 \\
Georgetown University & 1 \\
Louisiana State University & 1 \\
Massachusetts Institute of Technology & 1 \\
Mississippi State University & 1 \\
Northwestern University & 1 \\
Purdue University & 1 \\
Rutgers University & 1 \\
Stony Brook University & 1 \\
Saint Louis University & 1 \\
University of Arkansas & 1 \\
University of California, Los Angeles & 1 \\
University of Cincinnati & 1 \\
University of Connecticut & 1 \\
University of Florida & 1 \\
University of Glasgow, UK & 1 \\
University of Havana & 1 \\
University of Illinois Urbana--Champaign & 1 \\
University of Missouri & 1 \\
University of Nigeria, Nsukka & 1 \\
University of North Carolina at Charlotte & 1 \\
University of North Carolina at Chapel Hill & 1 \\
University of Oklahoma & 1 \\
University of Southern California & 1 \\
University of Texas at Austin & 1 \\
Virginia Polytechnic Institute and State University (Virginia Tech) & 1 \\
Wayne State University & 1 \\
Wesleyan University & 1 \\
Yale University & 1 \\
\end{longtable}

\vspace{0.5in}

\section{PhD Thesis Advisors}
This table lists the primary and secondary (if applicable) Ph.D. thesis advisors, as identified through the Mathematics Genealogy Project, MAD Pages, institutional records, and other publicly available records.

\begin{longtable}{p{0.78\textwidth}r}
\caption{PhD advisors of Black women mathematicians represented in the dataset.}
\label{tab:full_phd_advisors}\\
\toprule
\textbf{PhD Advisor} & \textbf{Number of Advisees} \\
\midrule
\endfirsthead

\multicolumn{2}{c}{\tablename\ \thetable\ -- continued from previous page}\\
\toprule
\textbf{PhD Advisor} & \textbf{Number of Advisees} \\
\midrule
\endhead

\midrule
\multicolumn{2}{r}{Continued on next page}\\
\endfoot

\bottomrule
\endlastfoot

Neil Hindman & 9 \\
James Dolan Reid & 4 \\
Alan Hopenwasser & 2 \\
Clinton G. McCrory & 2 \\
Doron Levy & 2 \\
Edward DeRoche & 2 \\
Frank Charles Hoppensteadt & 2 \\
Gilbert Baumslag & 2 \\
Mark William Brittenham & 2 \\
Mette Olufsen & 2 \\
Norman Levine & 2 \\
Richard Scott Pierce & 2 \\
Talmage James Reid & 2 \\
W E Deskins & 2 \\
William H. Kazez & 2 \\
David Constantine & 1 \\
Alex V. Kontorovich & 1 \\
Alexandra B Smirnova & 1 \\
Allan Peterson & 1 \\
Allison McCulloch & 1 \\
Alvin Bayliss & 1 \\
Aubrey Edward Landry & 1 \\
Benjamin Charles Brewster & 1 \\
Beth Herbel-Eisenmann & 1 \\
Bonnie Jean Dorr & 1 \\
Brian R. Hunt & 1 \\
Bruce Leonard Golden & 1 \\
Bruce P. Ayati & 1 \\
C. Einar (Carl) Hile & 1 \\
Calvin Cooper Moore & 1 \\
Chris Gennings & 1 \\
Christopher Kribs & 1 \\
Christopher Leininger & 1 \\
Colleen Catharine Mitchell & 1 \\
Curtis Claude McKnight & 1 \\
Daniel D. Anderson & 1 \\
David C. Wheeler & 1 \\
David W. Scott & 1 \\
Dennis Earl Davenport & 1 \\
Dianne O'Leary & 1 \\
Donald Beshers Parker & 1 \\
Donald Elton Edmondson & 1 \\
Douglas Grassel Dickson & 1 \\
Eduard Harabetian & 1 \\
Eric Saxton Brussel & 1 \\
Eric Victor Slud & 1 \\
Estelle Russek-Cohen & 1 \\
Eugenia Kalnay & 1 \\
Florence Sullivan & 1 \\
Frank Riedel & 1 \\
Frederick Russel Adler & 1 \\
Frederik S. Herzberg & 1 \\
Gabor Hetyei & 1 \\
Gary Christiansen & 1 \\
George Yuri Rabinovich Rainich & 1 \\
Herbert W. Hethcote & 1 \\
Hien T. Train & 1 \\
Hollylynne Lee & 1 \\
Howard Elman & 1 \\
Huixia Wang & 1 \\
Hyonho Chun & 1 \\
Isabel K. Darcy & 1 \\
Jack Elton Ohm & 1 \\
James A.M. Álvarez & 1 \\
James Francis Slegrade & 1 \\
Jianzhong Su & 1 \\
John A. Terilla & 1 \\
John Francis Monahan & 1 \\
John H. Hodges & 1 \\
John Joseph Benedetto & 1 \\
John Neuberger & 1 \\
John S. Mac Nerney. & 1 \\
Johny Lee Henderson & 1 \\
Jorge Martinez & 1 \\
Karen Hollebrands & 1 \\
Karin Saskia Dorman & 1 \\
Katherine Bennett Ensor & 1 \\
Kazufumi Ito & 1 \\
Kenneth Alexander Brown & 1 \\
Laura Jean Sheppardson & 1 \\
Leslie Hogben & 1 \\
Li Wang & 1 \\
Louis Shapiro & 1 \\
Luis Kramarz & 1 \\
Lynn Erbe & 1 \\
Maleafisha Joseph Pekwa Stephen Tladi & 1 \\
Manuel Rabina & 1 \\
Maria Korkina Cameron & 1 \\
Marius N. Nkashama & 1 \\
Martha Makowski & 1 \\
Martin Maxey & 1 \\
Marvin Knopp & 1 \\
Michael Owen Ball & 1 \\
Mohsen Razzaghi & 1 \\
Nozer Darabsha Singpurwalla & 1 \\
Oguz Durumeric & 1 \\
Olga Vitek & 1 \\
Óscar Chávez & 1 \\
Paola Sztajin & 1 \\
Paul John Smith & 1 \\
Paul Valentine Reichelderfer & 1 \\
Peter Louise Clark & 1 \\
Phillip Kutzko & 1 \\
Piotr Blass & 1 \\
Qingying Bu & 1 \\
Rainer Steinwandt & 1 \\
Ralph Conover Smith & 1 \\
Raymond Freese & 1 \\
Richard Mansfield Dudley & 1 \\
Robert John Sacker & 1 \\
Roger B. Sidje & 1 \\
Rudy Guerra & 1 \\
Russell Sherman Lehman & 1 \\
Ryan Robert Martin & 1 \\
Sam Taylor & 1 \\
Scott W Hansen & 1 \\
Shoshichi Kobayashi & 1 \\
Sivalingam Sritharan & 1 \\
Srikanth B Iyengar & 1 \\
Stanca Ciupe & 1 \\
Stephen Finch & 1 \\
Susan Marie Hermiller & 1 \\
Ted Belytschko & 1 \\
Thomas Roy Brahana & 1 \\
Thurlow Adrean Cook & 1 \\
Tong Li & 1 \\
Tony Fan Cheong Chan & 1 \\
Trevor Evans & 1 \\
Ulrich Krahmer & 1 \\
Warren Douglas Nichols & 1 \\
Wojciech Czaja & 1 \\
Xin Dang & 1 \\
Yong Wang & 1 \\
Yung Sze Choi & 1 \\
\end{longtable}

\noindent\textit{Note.} PhD thesis advisor information was unavailable for 9 women in the dataset.

\bibliography{biblio}


\end{document}